 \documentclass{amsart}
\usepackage{graphicx}       

\usepackage[latin1]{inputenc}
\usepackage[T1]{fontenc}
\usepackage{microtype} 

\usepackage{cite}
\usepackage{array,colortbl}
\usepackage{amsmath,amsfonts,amssymb,bm,amsthm,mathtools,mathrsfs}
\DeclareFontFamily{U}{mathx}{\hyphenchar\font45}
\DeclareFontShape{U}{mathx}{m}{n}{
      <5> <6> <7> <8> <9> <10>
      <10.95> <12> <14.4> <17.28> <20.74> <24.88>
      mathx10
      }{}
\DeclareSymbolFont{mathx}{U}{mathx}{m}{n}
\DeclareFontSubstitution{U}{mathx}{m}{n}
\DeclareMathAccent{\widecheck}{0}{mathx}{"71}

\def\citep#1#2{\cite[{#1}]{#2}}

\theoremstyle{plain}
 \newtheorem{theorem}{Theorem}
 \newtheorem{proposition}{Proposition}
 
 \newtheorem{corollary}{Corollary}
 
\theoremstyle{definition}
 \newtheorem{definition}{Definition}
 
\theoremstyle{remark}
 \newtheorem{remark}{Remark}

\usepackage{algorithm}
\usepackage{algorithmic}

\newcommand{\C}{{\mathbb{C}}} 
\newcommand{\N}{{\mathbb{N}}} 
\newcommand{\R}{{\mathbb{R}}} 
\newcommand{\Z}{{\mathbb{Z}}} 

\DeclareSymbolFont{bbold}{U}{bbold}{m}{n}
\DeclareSymbolFontAlphabet{\mathbbold}{bbold}

\begin{document}

\title{QMC strength for some random configurations on the sphere}

\author{V\'ictor de la Torre}
\address{Departament de Matem\`atiques i Inform\`atica, Universitat de Barcelona, Gran Via, 585, 08007 Barcelona, Spain}
\email{{\texttt{delatorre@ub.edu }}}
\author{Jordi Marzo}
\address{Departament de Matem\`atiques i Inform\`atica, Universitat de Barcelona, Gran Via, 585, 08007 Barcelona, Spain \vskip 0.01cm
Centre de Recerca Matem\`atica, Edifici C, Campus Bellaterra, 08193 Bellaterra, Spain}
\email{{\texttt{jmarzo@ub.edu}}}

\thanks{The authors have been partially supported by the Generalitat de Catalunya (grant 2017 SGR 359) and
the Spanish Ministerio de Ciencia, Innovaci\'on y Universidades (project MTM2017-83499-P)}

\begin{abstract}
A sequence $(X_N)\subset \mathbb{S}^d$ of $N$-point sets from the $d$-dimensional sphere has QMC strength $s^*>d/2$ if it has worst-case error of optimal order, $N^{s/d},$ for Sobolev spaces of order $s$ for all $d/2<s<s^*,$ and 
the order is not optimal for $s> s^*.$
In \cite{BrSaSlWo} conjectured values of the strength are given for some well known point families in $\mathbb S^2$ based on numerical results. We study 
the average QMC strength for some related random configurations.
\end{abstract} 

\maketitle

\section{Introduction and main results}
Let $\mathbb{S}^{d}=\{x\in \R^{d+1}\; : \; |x|=1\}$ be the unit sphere with the normalized Lebesgue measure $\sigma_d.$ Given an integer $\ell\ge 0$, let $\mathcal{H}_\ell$ be the vector space of 
the spherical harmonics of degree $\ell,$ that is, the space of eigenfunctions of the Laplace-Beltrami operator, $\Delta,$ with eigenvalue $\lambda_\ell=\ell (\ell +d-1)$,
\begin{equation*}
-\Delta Y=\lambda_\ell Y,\;\;\;\; Y\in \mathcal{H}_{\ell}.
\end{equation*}
We denote $h_\ell=\mbox{dim}\, \mathcal{H}_\ell\sim \ell^{d-1}.$ The space of spherical harmonics of degree at most 
$L$ in $\mathbb{S}^d$ is
$$\mathcal{P}_L(\mathbb{S}^d)=\bigoplus_{\ell=0}^L \mathcal{H}_\ell,$$
which is also the space of polynomials of degree at most $L$ in 
$\R^{d+1}$ restricted to $\mathbb{S}^d.$  Clearly $d(L):=\mbox{dim}\, \mathcal{P}_L(\mathbb{S}^d) \sim L^d.$

Let $L^2(\mathbb{S}^d)$ be the Hilbert space of square integrable real functions in $\mathbb{S}^d$
with the inner product
\[
\langle f,g \rangle=\int_{\mathbb{S}^d} f(x)g(x)\,d\sigma_d(x),\;\;\; f,g\in L^2(\mathbb{S}^d).
\]
One has that $L^2(\mathbb{S}^d)=\bigoplus_{\ell\ge 0} \mathcal{H}_\ell$ and the Fourier series expansion of a function $f \in L^2(\mathbb{S}^d)$ is given by
\begin{equation*}
f=\sum_{\ell,k} f_{\ell,k} Y_{\ell,k},\qquad f_{\ell,k}=\langle f,Y_{\ell,k} \rangle=\int_{\mathbb{S}^d} f\, Y_{\ell,k} \,d\sigma,
\end{equation*}
where $\{ Y_{\ell,k} \}_{k=1}^{h_\ell}$ is an orthonormal basis of $\mathcal{H}_\ell.$ 

Following \cite{BrSaSlWo} we define, for $s\ge 0,$ the $L^2(\mathbb{S}^d)$-based Sobolev spaces of order $s$ as the Hilbert space
$$\mathbb{H}^{s}(\mathbb{S}^d)=\left\{ f\in L^2(\mathbb{S}^d)\;\;:\;\; 
\sum_{\ell=0}^{+\infty}\sum_{k=1}^{h_\ell} (1+\lambda_\ell)^s |f_{\ell,k}|^2<+\infty \right\},$$
with the norm
$$\| f \|_{\mathbb{H}^{s}(\mathbb{S}^d)}=\left( \sum_{\ell=0}^{+\infty}\sum_{k=1}^{h_\ell} \frac{1}{a_\ell^{(s)}} |f_{\ell,k}|^2  \right)^{1/2},$$ 
where $a_\ell^{(s)}\sim (1+\ell^2)^{-s}.$ It is well known that 
$\mathbb{H}^{s}(\mathbb{S}^d)$ is continuously embedded in $\mathcal{C}(\mathbb{S}^d)$ if $ s >d/2$ and it has, in this range, a reproducing kernel given by
$$K^{(s)}(x,y)=K^{(s)}(x\cdot y)=\sum_{\ell=0}^\infty \sum_{k=0}^{h_\ell} a_\ell^{(s)} Y_{\ell,k}(x) Y_{\ell,k}(y),$$
i.e. for all $x\in \mathbb{S}^d$ and $f\in \mathbb{H}^{s}(\mathbb{S}^d)$,
$$f(x)=\langle f, K^{(s)}(x,\cdot) \rangle_{\mathbb{H}^{s}(\mathbb{S}^d)}.$$

 \begin{definition}
Given $s>d/2,$ a sequence $(X_N)$ of $N-$point configurations $X_N\subset \mathbb{S}^d$ is a sequence of QMC designs for $\mathbb{H}^{s}(\mathbb{S}^d)$ (or an $s$-QMC design) if there exists $C_{d,s}>0$ such that for all $N\ge 1$
\begin{equation}\label{DeMa_wcedef}
\mbox{wce}(X_N,\mathbb{H}^{s}(\mathbb{S}^d))\le C_{d,s} N^{-s/d},
\end{equation}
where the worst-case error 
is defined by
$$\mbox{wce}(X_N,\mathbb{H}^{s}(\mathbb{S}^d))=
\sup_{\| f \|_{\mathbb{H}^{s}(\mathbb{S}^d)}\le 1} \left\{ 
\left| \frac{1}{N}\sum_{x\in X_N}f(x)-\int_{\mathbb{S}^d}f(x)d\sigma(x) \right| \; : \; f\in \mathbb{H}^{s}(\mathbb{S}^d) \right\}.$$
 \end{definition}

From \cite[Theorem 3.1]{BrChCoGiSeTr} it follows that if $(X_N)$ is a sequence of QMC designs for $\mathbb{H}^s(\mathbb{S}^d)$ it is also a QMC design for all $\mathbb{H}^{s'}(\mathbb{S}^d)$ for $\frac{d}{2}<s'<s,$ see also \cite[Lemma 23]{BrSaSlWo}.
The maximal $s^*>\frac{d}{2}$ where (\ref{DeMa_wcedef}) holds for all $s<s^*$ is the QMC strength of the sequence $(X_N),$ \cite{BrSaSlWo}. The strength can be seen as a measure of the regularity of the sequence.
In the definition above $(X_N)$ may be defined only for a subsequence of the natural numbers $N$ converging to $+\infty.$  
 
The problem of determining the strength for a given sequence seems to be quite difficult and its value has been determined only in a few cases. 
It was shown by Hesse and Sloan \cite{HeSl} and by Brauchart and Hesse \cite{BrHe} that sequences of optimal quadrature formulas satisfying some regularity property, true in particular for quadrature formulas with positive weights, are
$s$-QMC designs for all $s>d/2$, i.e., 
they have $s^*=+\infty.$ It was observed in \cite{BrSaSlWo} that this previous result, together with the existence of optimal spherical designs \cite{BoRaVi}, implies the existence of spherical designs with strength $s^*=+\infty.$
Also in \cite[Theorem 14]{BrSaSlWo}, the authors show that maximizers of the sum of suitable powers of the Euclidean distance between pairs of points are $s$-QMC designs for all $d/2<s<d/2+1.$ To the best of our knowledge these are the only 
cases where the strength is known.

Values for the strength where conjectured in \cite{BrSaSlWo} for some well known point configurations in $\mathbb S^2$ based on numerical results. In particular, for Fekete points the conjecture is $s^*=3/2,$ for equal 
area points $s^*=2$ and for minimal logarithmic energy points $s^*=3,$  see next section for definitions. The expected worst-case error of some random configurations was also studied in \cite{BrSaSlWo}.
 \begin{definition}
Let $(X_N)$ be a sequence of random $N-$point configurations on $\mathbb{S}^d$ following some distribution and let $s>d/2.$ As in the deterministic case, we allow the subindex $N$ to follow a subsequence converging to $+\infty.$ We say 
that $(X_N)$ is a sequence of  QMC designs for $\mathbb{H}^{s}(\mathbb{S}^d)$ on average (or $s$-QMC designs on average) if there 
exists $C_{d,s}>0$ such that
\begin{equation}\label{DeMa_wcerandef}
\sqrt{\mathbb E [\mbox{wce}(X_N,\mathbb{H}^{s}(\mathbb{S}^d))^2]}\le C_{d,s} N^{-s/d}.
 \end{equation}
 \end{definition}

 With obvious changes in the proof of Lemma 23 in \cite{BrSaSlWo}, the following property shows that if a sequence of random point configurations is a sequence of $s$-QMC designs on average, for some $s>d/2,$ then it is also 
 of $s'$-QMC designs for $d/2<s'<s.$
\begin{proposition}\label{DeMa_lowranQMC}
 Given $s>d/2$, if $\mathbb E [\mbox{\rm wce}(X_N,\mathbb{H}^{s}(\mathbb{S}^d))^2]\le 1$, then there exists a constant $C_{d,s',s}>0$ such that
 $$\mathbb E [\mbox{\rm wce}(X_N,\mathbb{H}^{s'}(\mathbb{S}^d))^2]\le C_{d,s',s} \left(\mathbb E [\mbox{\rm wce}(X_N,\mathbb{H}^{s}(\mathbb{S}^d))^2]\right)^{s'/s}, \;\;\frac{d}{2}<s'<s.$$
\end{proposition}
This last result allows us to define the average QMC strength, as in the deterministic case, as the maximal value of $s> d/2$ for which (\ref{DeMa_wcerandef}) holds.

It was shown in \cite[Theorem 7]{BrSaSlWo} that uniform i.i.d. points on the sphere are not $s$-QMC designs on average for any $s>d/2.$ On the other hand, 
\cite[Theorem 21,Theorem 22]{BrSaSlWo} shows that points from jittered 
sampling (i.e. uniform i.i.d point taken with respect to an area regular partition) have average strength 
$d/2+1.$ Observe that in this last case the average strength matches the numerically conjectured in \cite{BrSaSlWo} for the related equal area points in $\mathbb{S}^2.$ 

%


\subsection{Harmonic ensemble}

In our first result we show that points from the harmonic ensemble have average strength $\frac{d+1}{2}.$ The mode of this distribution corresponds to the Fekete points, see next section, for which it was conjectured
strength $3/2$ in \cite{BrSaSlWo}. The expected worst case error of this process was previously studied in \cite{Hi}.

\begin{theorem}\label{DeMa_wce_harmonic}
Let $(X_N)$ be a sequence where $X_N$ is an $N-$point set drawn from the harmonic ensemble in $\mathbb{S}^d.$ Observe that $N$ must be of the form $d(L)$ for some natural $L.$ Then $(X_N)$
forms a sequence of QMC designs on average for $\frac{d}{2}<s<\frac{d+1}{2}.$ Moreover
\begin{equation}\label{DeMa_divergence_harmonic}
\lim_{N\to +\infty }N^{\frac{d+1}{d}}\mathbb E [ \mbox{\rm wce} (X_N;\mathbb{H}^{\frac{d+1}{2}}(\mathbb{S}^2))^2 ]=+\infty, 
\end{equation}
therefore $(X_N)$ is not a QMC design on average if $s>\frac{d+1}{2}$ and the average QMC strength is $\frac{d+1}{2}.$
\end{theorem}

For the harmonic ensemble we can deduce, from results in \cite{BrSaSlWo,BeMaOr}, see also \cite{Be}, almost sure optimality of the worst-case error up to a logarithmic
factor.

\begin{corollary} \label{DeMa_cor_as}
  For every $M>0$ and $\frac{d}{2}<s<\frac{d+1}{2}$, there exists $C_{d,s,M}>0$ such that
  \begin{equation}\label{DeMa_bounds_prob}
  \mathbb P\left( \mbox{\rm wce} (X_N;\mathbb{H}^{s}(\mathbb{S}^d)) \le C_{d,s,M} \frac{(\log N)^{\frac{2s}{d+1}}}{N^{\frac{s}{d}}}  \right)\ge 1-\frac{1}{N^M}, 
  \end{equation}
where $X_N$ is an $N$-point set drawn from the harmonic ensemble. Therefore, for fixed $\frac{d}{2}<s<\frac{d+1}{2}$ there exists $C_{d,s}>0$ such that, with probability 1 and for $N$ large enough,
$$\mbox{\rm wce} (X_N;\mathbb{H}^{s}(\mathbb{S}^d)) \le C_{d,s} \frac{(\log N)^{\frac{2s}{d+1}}}{N^{\frac{s}{d}}}.$$
\end{corollary}


\subsection{Spherical ensemble}

The spherical ensemble is also a determinantal point process on the sphere $\mathbb{S}^2.$ Applying results from    
\cite{AlZa} it was shown in \cite{Hi} that points from the spherical ensemble are $s$-QMC designs on average for $1<s<2.$ One can easily see that 2 is indeed the average strength. Observe that 
the mode of this distribution is the set of elliptic Fekete points, i.e. minimizers of the logarithmic energy. In this case there is no coincidence with the conjectured strength from \cite{BrSaSlWo}, which was $3.$

\begin{theorem}\label{DeMa_wce_spherical}
Let $(X_N)$ be a sequence where $X_N$ is an $N-$point set drawn from the spherical ensemble. Then $(X_N)$ is a sequence of QMC designs on average for $1<s<2,$ and for $s\in (2,3)$
there exists a constant $C>0$ such that
\[N^2\mathbb{E}\left[\mbox{\rm wce}(X_N;\mathbb{H}^s(\mathbb{S}^2))^2\right]\geq C, \]
i.e. the average strength is $s^{*}=2$.
\end{theorem}


\subsection{Zeros of Elliptic polynomials}

In our last result, by using the expected energy expansions from \cite{DeMa} we prove that the average strength for the zeros of the elliptic polynomials behaves better than for all these previous random processes and coincides with the 
conjectured strength in \cite{BrSaSlWo}
for the logarithmic energy minimizers. 

\begin{theorem}\label{DeMa_wce_theorem}
Let $(X_N)$ be a sequence where $X_N$ is an $N-$point set drawn from zeros of elliptic polynomials mapped to the sphere by the stereographic projection. Then $(X_N)$ is a sequence of $s$-QMC designs on average for $1<s<3,$ and for $s\in (3,4)$
there exists a constant $C>0$ such that
\[N^3\mathbb{E}\left[\mbox{\rm wce}(X_N;\mathbb{H}^s(\mathbb{S}^2))^2\right]\geq C, \]
i.e. the average strength is $s^{*}=3$.
\end{theorem}

\begin{remark}
A result similar to Corollary \ref{DeMa_cor_as} can be proved also for configurations given by determinantal point processes, like the jittered sampling \cite{BrGrKuZi} and the spherical ensemble \cite{AlZa}, using the 
concentration results for determinantal point processes in \cite{PePe}. The bounds are far from sharp. 
For the spherical ensemble an almost optimal bound has been proved using a concentration of measure inequality particular of the spherical ensemble, \cite{Be}.
\end{remark}

\section{Background}


\subsection{Riesz energy and worst-case error}
The Riesz or logarithmic energy of a set of $N$ different points $x_1,\ldots,x_N$
on the unit sphere $\mathbb{S}^d \subset \mathbb R^{d+1}$ is 
$$E_s(x_1,\dots ,x_N)=\sum_{i\neq j}f_s(|x_i-x_j|),$$ 
where $f_s(r)=r^{-s}$ for $s\neq 0$ and $f_0(r)=-\log r$ are, respectively, the Riesz and logarithmic potentials.
From now on, to simplify the notation, we write $E_s$ for $E_s(x_1,\dots ,x_N)$ when the set of points is clear from the context.
This quantity has a continuous version for measures which for the normalized surface measure $\sigma_d$ and  $0\neq s < d$ is
$$V_s(\mathbb{S}^d)=\int_{\mathbb{S}^d}\int_{\mathbb{S}^d}f_s(|x-y|)\,d\sigma(x)d\sigma(y)=2^{d-s-1}\frac{\Gamma\left(\frac{d+1}{2}\right)\Gamma\left(\frac{d-s}{2}\right)}{\sqrt{\pi}\,\Gamma\left(d-\frac{s}{2}\right)}.$$
In \cite{BrSaSlWo}, the authors obtained a formula for the worst-case error of an $N$-point set in terms of Riesz energies, provided $s$ is not a positive integer. When $d/2<s<d/2+1,$ the expression reads
\begin{equation}\label{DeMa_s12}
\mbox{wce}(X_N,\mathbb{H}^{s}(\mathbb{S}^d))^2=-\frac{1}{N^2}\left(E_{d-2s}-V_{d-2s}(\mathbb{S}^d)N^2\right),
\end{equation}
while for $d/2+M<s<d/2+M+1$, with $M$ a positive integer,
\begin{equation}\label{DeMa_others}
\mbox{wce}(X_N,\mathbb{H}^{s}(\mathbb{S}^d))^2=\frac{1}{N^2} \left[\sum_{j,i=1}^N \mathcal{Q}_M(x_j \cdot x_i)+(-1)^{M+1}(E_{d-2s}-V_{d-2s}(\mathbb{S}^d)N^2)\right],
\end{equation}
where
\[\mathcal{Q}_M( x_j \cdot x_i):=\sum_{l=1}^M\left((-1)^{M+1-l}-1\right)\alpha_l^{(s)}h_l P_l^{(d)}( x_j \cdot x_i),\]
with $P_l^{(d)}(x)$ the Gegenbauer polynomial normalized by $P_l^{(d)}(1)=1$, 
\begin{equation}\label{DeMa_defalpha}
\alpha_l^{(s)}=V_{d-2s}(\mathbb{S}^d)\frac{(-1)^{M+1}(1-s)_l}{(1+s)_l},
\end{equation}
and $(\cdot)_n$ for $n\in \N$ is the Pochhammer symbol given, for $x\notin \Z_{\leq 0},$ by $(x)_0=1$ and 
$$(x)_n=\frac{\Gamma(x+n)}{\Gamma(x)},\;\;n\ge 1.$$

A way to measure the uniformity of a finite set of points $X_N=\{x_1, \dots, x_N\}\subset \mathbb{S}^d$ is to consider discrepancies. In particular, we are going to consider the $L^\infty$ spherical cap discrepancy
\[ 
\mathbb{D}_\infty(X_N)=\sup_{x\in \mathbb{S}^d,r>0} \Bigl| \frac{1}{N}\sum_{i=1}^{N} \chi_{D(x,r)}(x_i)-\sigma_d(D(x,r)) \Bigr|,
\]
and the $L^2$ averaged version
$$\mathbb{D}_2(X_N)=\left( \int_0^\pi \int_{\mathbb{S}^d} \left|   \frac{1}{N}\sum_{i=1}^{N} \chi_{D(x,r)}(x_i)-\sigma_d(D(x,r))   \right|^2 d\sigma_d(x)\sin r dr  \right)^{1/2}.$$
It is well known that a sequence of point sets $(X_{N_k})_k$ is asymptotically uniformly distributed 
if and only if $\lim_{k\to +\infty} \mathbb{D}_p (X_{N_k})=0,$ for $p=\infty$ or $p=2,$ \cite[Section 6.1]{BoHaSa}.

Observe that it follows from (\ref{DeMa_s12}) that for $s=(d+1)/2$ there is a direct relation between energy, discrepancy and worst-case error 
\begin{align}\label{DeMa_equality_extended}
\mbox{\rm wce}  (X_N;\mathbb{H}^\frac{d+1}{2}(\mathbb{S}^d)) &=\left(V_{-1}(\mathbb{S}^d)-\frac{1}{N^2}\sum_{i,j=0}^N |x_i-x_j|\right)^{1/2} \nonumber
\\
&
=
\sqrt{\frac{d\sqrt{\pi}\Gamma(\frac{d}{2})}{\Gamma(\frac{d+1}{2})}} \mathbb{D}_2(X_N), 
\end{align}
where the last equality follows from Stolarsky's invariance principle, see \cite{BrDi,BrSaSlWo}.

\subsection{Point processes}

A point process $\mathcal{X}$ on the sphere $\mathbb S^d$
is a random integer-valued positive Radon measure on $\mathbb{S}^d$ that almost surely assigns at most measure 1 to singletons, for general definitions and more backgound see \cite{HoKrPeVi,AnGuZe}.
A random point process in $\mathbb{S}^d$ can be identified with a random discrete subset of $\mathbb{S}^d.$ A way to 
describe a random point process is to specify the random variable counting the number of points of the process in $D$,
for all Borel sets $D\subset \mathbb{S}^d$. We denote this random variable as $\mathcal{X}(D).$ For all the point processes we consider, the distribution is characterized by the so-called joint intensity functions. The joint 
intensities with respect to some background measure $\mu$ in $\mathbb S^d$ are functions $\rho_k(x_1,\dots, x_k)$
that for any family of mutually disjoint subsets $D_1,\dots ,D_k \subset \mathbb{S}^d$ satisfy
\[
\mathbb{E}\left[ \mathcal{X}(D_1)\cdots \mathcal{X}(D_k) \right] 
=\int_{D_1\times \dots \times D_k} \rho_k(x_1,\dots, x_k) d\mu(x_1)\dots 
d\mu    (x_k).
\]
We assume that $\rho_k(x_1,\dots, x_k)=0$ when $x_i=x_j$ for $i\neq j$. Note that if the distribution of the process is given by the joint densities $p(x_1,	\dots , x_N)$ then
$$\rho_k(x_1,\dots, x_k)=\frac{N!}{(N-k)!}\int_{\R^{N-k}}p(x_1,\dots, x_N)dx_{k+1}\dots dx_{N},$$
see \cite[p. 10]{HoKrPeVi}.

A particular type of random point process are the so-called determinantal point processes (DPP) introduced by Macchi in 1975, \cite{Ma}. 
A random point process on the sphere is called determinantal with kernel 
$K:\mathbb{S}^d\times \mathbb{S}^d \rightarrow \C$ if the joint intensities
with respect to a background measure $\mu$ are given by
\[
\rho_k(x_1,\dots, x_k)=\det (K(x_i,x_j))_{1\le i,j\le k}, 
\]
for every $k\ge 1$ and $x_1,\dots, x_k \in \mathbb{S}^d$. For the determinantal point process to exist it is enough that the 
kernel function
$K:\mathbb{S}^d\times \mathbb{S}^d \rightarrow \C$ determines a self-adjoint integral operator in $L^2(\mathbb S^d)$ that is locally of trace class, \cite[Theorem 4.5.5]{HoKrPeVi}.



\subsection{Harmonic ensemble} The harmonic ensemble is the determinantal point process in $\mathbb{S}^d$ with $d(L)$ points a.s.
induced by the reproducing kernel of the space $\mathcal P_L(\mathbb S^d)$  of polynomials in $\R^{d+1}$ of degree at most $L$, 
\[
K_L(x,y)=\frac{d(L)}{\binom{L+\frac{d}{2}}{L}}P_L^{(1+\lambda,\lambda)}(\langle x,y 
\rangle),\;\;x,y\in \mathbb{S}^d,  
 \]
where $\lambda=\frac{d-2}{2}$ and the Jacobi polynomials 
$P_L^{(1+\lambda,\lambda)}(t)$ are normalized to have $P_L^{(1+\lambda,\lambda)}(1)=\binom{L+\frac{d}{2}}{L}.$

Observe that the mode of the distribution, in some sense the value that appears most often in a set of data values sampled from this DPP, is the maximum of the joint density 
$$p(x_1,\dots ,x_{d(L)})=(1/d(L)!) \det(K_L(x_i,x_j))_{1\le i,j\le d(L)},$$
for $x_1,\dots ,x_{d(L)}\in \mathbb S^d.$ The mode points are therefore the Fekete points that maximize the absolute value of the determinant 
$$V_L(x_1,\dots , x_{d(L)})=\left|\det(\phi_i(x_j))_{1\le i,j\le d(L)}\right|,$$
for $\phi_1,\dots, \phi_{d(L)}$ a basis of the space $\mathcal{P}_L(\mathbb S^d).$ Fekete points are called also extremal fundamental systems, \cite{Re}.  Sloan and Womersley conjectured that they have all positive 
cubature weights \cite{SlWo} and were shown to be asymptotically uniformly equidistributed in \cite{MaOr,BeBoNy}. 


\subsection{Spherical ensemble}

Given $A,B$ independent $n\times n$ random matrices with i.i.d. complex standard entries, the eigenvalues $\lambda_1,\dots ,\lambda_n$ of $A^{-1}B$ form a DPP in $\C$ with kernel $(1+z\overline{w})^{n-1}$ with respect to the 
background measure $\frac{n}{\pi(1+|z|^2)^{n+1}}dm(z),$ \cite{Kr}. When mapped to the sphere $\mathbb S^2$ through the stereographic projection $g(x_i)=\lambda_i$ the points  have density 
$$C \prod_{i\neq j}|x_i-x_j|^2,$$
with respect to the surface measure in $\mathbb S^2.$ This shows that the mode of this point process on the sphere is given by the minimizers of the logarithmic energy
$$E_0(x_1,\dots, x_N)=\sum_{i\neq j} \log\frac{1}{|x_i-x_j|},$$
i.e. for the (elliptic) Fekete points, \cite{BoHaSa}.

In \cite{AlZa} the authors obtain, among other results, the expected Riesz and 
logarithmic energies for points of this process. In particular, they prove that for $x_1,\dots, x_N\in \mathbb S^2$ points drawn from the spherical ensemble, if $s<4$, $s\neq 0,2$,
\begin{equation}\label{DeMa_energy_sph}
\mathbb{E}\left[E_s\right]=\frac{2^{1-s}}{2-s}N^2-\frac{\Gamma(N)\Gamma(1-s/2)}{2^s\Gamma(N+1-s/2)}N^2.
\end{equation}

\subsection{Zeros of Elliptic polynomials} The last point process that we consider is given by the images in $\mathbb S^2$ by the stereographic 
projection of the zeros  of the elliptic polynomials (also called Kostlan-Shub-Smale or $SU(2)$ polynomials)
\begin{equation*}
\sum_{n=0}^N a_n\sqrt{\binom{N}{n}}z^n,
\end{equation*}
where $a_n$ are i.i.d. random variables with standard complex Gaussian distribution. The study of the zeros of these polynomials started in quantum chaotic dynamics \cite{BeBoLe92,BeBoLe96}. Because of the 
connection with minimal logarithmic energy points on the sphere, unveiled in \cite{SmSh}, it was natural to study the expected logarithmic energy. This was done in \cite{ArBeSh} where the authors
obtained the following closed expression for the expected logarithmic energy
\begin{equation} \label{DeMa_armentano}
\mathbb E [ E_0 ]=\left(\frac{1}{2}-\log 2\right)N^2
-\frac{1}{2}N\log N
-\left(\frac{1}{2}-\log 2\right)N.
\end{equation} 
The extension of this result to the Riesz $s$-energy for $s<4$ has been recently obtained in \cite{DeMa}. In order to study the expected QMC strength we will use that
for $x_1,\dots, x_N\in \mathbb S^2$ points drawn from zeros of elliptic polynomials mapped to the sphere by the stereographic projection and $s<4$, $s\neq 0,2$,
\begin{equation}\label{DeMa_weakasymptotics}
\mathbb{E}[E_s]=\frac{2^{1-s}}{2-s}N^2+C(s)N^{1+s/2}+o_{N\to\infty}(N^{1+s/2}),
\end{equation}
where
\[C(s)=\frac{1}{2^s}\frac{s}{2}\left(1+\frac{s}{2}\right)\Gamma\left(1-\frac{s}{2}\right)\zeta\left(1-\frac{s}{2}\right).\]
In the particular case $s=-2$,
\begin{equation}\label{DeMa_cas-2}
\mathbb{E}[E_{-2}]=2N^2-8\frac{\zeta(3)}{N}+o_{N\to\infty}\left(\frac{1}{N}\right).
\end{equation}

\section{Proofs }

\subsection{The harmonic ensemble (Theorem \ref{DeMa_wce_harmonic})} We start with some preliminaries in order to prove Theorem \ref{DeMa_wce_harmonic}.

Let $\frac{d}{2}<s<\frac{d}{2}+1$. From formula \eqref{DeMa_s12}, for 
$N=d(L)\sim L^d$ points drawn from the harmonic ensemble on the sphere we get
$$\mathbb{E}[\mbox{\rm wce} (X_N;\mathbb{H}^s(\mathbb{S}^d))^2]=\frac{1}{N^2}\int_{\mathbb{S}^d}\int_{\mathbb{S}^d}K_L(x,y)^2 |x-y|^{2s-d}d\sigma(x)d\sigma(y)$$
$$=\frac{d(L)^2}{P_L^{(1+\lambda,\lambda)}(1)^2 N^2}\int_{\mathbb{S}^d}\int_{\mathbb{S}^d} P_L^{(1+\lambda,\lambda)}(\langle x,y \rangle)^2|x-y|^{2s-d}d\sigma(x)d\sigma(y)$$
$$=\frac{C_d}{P_L^{(1+\lambda,\lambda)}(1)^2 } \int_{\mathbb{S}^d} P_L^{(1+\lambda,\lambda)}(\langle x,{\bf n} \rangle)^2|x-{\bf n}|^{2s-d}d\sigma(x)$$
$$=\frac{C_d}{ P_L^{(1+\lambda,\lambda)}(1)^2}
\int_{-1}^1 P_L^{(1+\lambda,\lambda)}(t)^2 (1-t)^{s-1} (1+t)^{\frac{d}{2}-1}dt,$$
where $\bf n$ stands for the north pole of $\mathbb{S}^d.$

From the asymptotic property of the gamma function 
\[
\lim_{n\to \infty} \frac{\Gamma(n+\alpha)}{\Gamma(n)n^\alpha}=1,\;\;\alpha\in 
\R,
\]
we get that
$$P_L^{(1+\lambda,\lambda)}(1)=\binom{L+\frac{d}{2}}{L}\sim \frac{1}{\Gamma(\frac{d}{2}+1)}L^{d/2}.$$
Therefore we have that for some constant $C_{d,s}>0$,
\begin{equation}\label{DeMa_closedform}
 \mathbb{E}[\mbox{\rm wce} (X_N;\mathbb{H}^s(\mathbb{S}^d))^2]=\frac{C_{d,s} }{L^d}
\int_{-1}^1 P_L^{(1+\lambda,\lambda)}(t)^2 (1-t)^{s-1} (1+t)^{\frac{d}{2}-1}dt.
\end{equation}

The following lemma is an extension to $-1<a<d$ of a result proved in \cite{BeMaOr} for $0<a<d.$


\begin{proposition}\label{DeMa_pr:integralestimate}
Given $-1<a<d$, 
\begin{equation*}
\lim_{L\to \infty} \frac{1}{L^a}\int_{-1}^1 P_L^{(1+\lambda,\lambda)}(t)^2 
(1-t)^{\lambda-\frac{a}{2}}(1+t)^{\lambda}\,dt=
2^{\frac{a}{2}+d}\int_0^{\infty} \frac{J_{1+\lambda}(t)^2}{t^{1+a}} dt
\end{equation*}
and the last integral converges.
\end{proposition}

\proof  
The proof is essentially the same as in \cite[Proposition 6]{BeMaOr} but with a few changes in the last estimates. We split the integral
\begin{align*}
 \int_{-1}^1 & L^{-a}  P_L^{(1+\lambda,\lambda)}(t)^2 
(1-t)^{\lambda-\frac{a}{2}}(1+t)^{\lambda}\,dt
\\
=&\left[\int_{-1}^{-\cos \frac{c}{L}}+\int_{-\cos \frac{c}{L}}^{\cos 
\frac{c}{L}}+
\int_{\cos \frac{c}{L}}^1 \right] L^{-a} P_L^{(1+\lambda,\lambda)}(t)^2 
(1-t)^{\lambda-\frac{a}{2}}(1+t)^{\lambda}\,dt
\\
=& A(c,L)+B(c,L)+C(c,L),
\end{align*}
where $c>0$ is fixed and $c<\pi L$. For the boundary parts we do a change of 
variables $t=\cos(x/L)$ to get
\begin{align*}
& C (c,L)  
\\
&
=2^{a/2}\int_0^c L^{-2-2\lambda}  P_L^{(1+\lambda,\lambda)}\left(\cos 
\frac{x}{L}\right)^2 \left( \frac{\sin \frac{x}{L}}{\frac{x}{L}} 
\right)^{2\lambda+1}
\left( \frac{1-\cos \frac{x}{L}}{\frac{1}{2}\left(\frac{x}{L}\right)^2} 
\right)^{-a/2} x^{2\lambda+1-a}\,dx.
\end{align*}
Using the Mehler-Heine asymptotic formula \cite[p. 192]{Sz} and the elementary limits
\[
 \lim_{L\rightarrow\infty}\frac{\sin \frac{x}{L}}{\frac{x}{L}}=1,\quad 
\lim_{L\rightarrow\infty}\frac{1-\cos 
\frac{x}{L}}{\frac{1}{2}\left(\frac{x}{L}\right)^2}=1,
\]
we conclude:
\[
 \lim_{L\rightarrow\infty}C(c,L)=2^{\frac{a}{2}+d}\int_0^c 
\frac{J_{1+\lambda}(x)^2}{x^{1+a}} dx.
\]

For the other end of the interval, using the change of variables $t=-\cos(x/L)$ 
we get
\[
A(c,L)=
\int_0^c L^{-2-2\lambda}  P_L^{(1+\lambda,\lambda)}\left(-\cos 
\frac{x}{L}\right)^2 
\left( \frac{\sin \frac{x}{L}}{\frac{x}{L}} \right)^{2\lambda+1}
\left(\frac{ 1+\cos \frac{x}{L} }{\left(\frac{x}{L}\right)^2}\right)^{-a/2} 
x^{2\lambda+1} dx,
\]
and using Mehler-Heine again this expression converges to zero when $L\to \infty$.
For the middle term we use classical asymptotic estimates of the Jacobi polynomials \cite[Theorem 8.21.13]{Sz} 
\begin{equation*}         			
    P_{L}^{(1+\lambda,\lambda)}(\cos
    \theta)=\frac{k(\theta)}{\sqrt{L}}\left\{ \cos \left((L+\lambda+1)\theta+\gamma 
\right)+\frac{O(1)}{L\sin
    \theta}\right\},
\end{equation*}
    if $c/L\le \theta \le \pi-(c/L)$ 
\[
k(\theta)=\pi^{-1/2} \left( \sin \frac{\theta}{2}\right)^{-\lambda-3/2}\left(\cos
    \frac{\theta}{2}\right)^{-\lambda-1/2},\quad \mbox{and}
\;\;\gamma=-\left(\lambda+\frac{3}{2}\right)\frac{\pi}{2}.
\]
We get
\begin{align*}
0\leq & B(c,L)  \lesssim \frac{1}{L^{a+1}}\int_{\frac{c}{L}}^{\pi-\frac{c}{L}} 
\frac{1}{(\sin \frac{\theta}{2})^{a+2}}d\theta\lesssim
\frac{1}{L^{a+1}}\int_{\frac{c}{L}}^{\pi-\frac{c}{L}} 
\frac{1}{\theta^{a+2}}d\theta
\\
&
\le 
\frac{1}{L^{a+1}}\int_{\frac{c}{L}}^{\pi} 
\frac{1}{\theta^{a+2}}d\theta=\frac{1}{L^{a+1}(a+1)}\left[ \frac{L^{a+1}}{c^{a+1}}-\frac{1}{\pi^{a+1}} \right]=
\frac{1}{c^{a+1}(a+1)}+o(L).
\end{align*} 

Finally, observe that close to zero $J_{1+\lambda}(x)\sim x^{1+\lambda}$ and 
$J_{1+\lambda}(x)\lesssim x^{-1/2}$ for big $x,$  so the integral above converges precisely for $-1<a<d.$

\qed



\proof(Theorem \ref{DeMa_wce_harmonic})
Now for any $\frac{d}{2}<s<\frac{d+1}{2}$ we apply the proposition above to $-1<t=d-2s<0$ and we get from (\ref{DeMa_closedform}) that 
there exist a (different) constant $C_{d,s}>0$ such that 
$$\lim_{L\to +\infty}L^{2s} \mathbb{E}[\mbox{\rm wce} (X_N;\mathbb{H}^s(\mathbb{S}^d))^2]=C_{d,s}.$$
This shows that points drawn from the harmonic ensemble form an $s$-QMC design on average
for $\frac{d}{2}<s<\frac{d+1}{2}.$

To get (\ref{DeMa_divergence_harmonic}) we use again (\ref{DeMa_closedform}) for $s=\frac{d+1}{2}$ and the representation of the integral in terms of generalized hypergeometric function \cite[p. 288]{ErMaObTr}
\begin{align*}
 \mathbb{E} & [\mbox{\rm wce} (X_N;\mathbb{H}^s(\mathbb{S}^d))^2]\sim L^{-d} \frac{\Gamma\left(\frac{1}{2}+L \right) \Gamma\left(\frac{d}{2}+L \right)\Gamma\left(\frac{d}{2}+L+1 \right)}{\Gamma(L+1)^2 \Gamma\left(d+\frac{1}{2}+L \right)^2}
 \\
 &
 \times
 {}_{4}F_{3}\left( 
-L,d+L,\frac{d+1}{2},\frac{1}{2};\frac{d}{2}+1,d+\frac{1}{2}+L,-L+\frac{1}{2} ;1  \right),
\end{align*}
where all the constants depend only on $d.$

It is easy to see (by induction) that the quotient
$$\frac{(-L)_n (d+L)_n}{(d+\frac{1}{2}+L)_n (-L+\frac{1}{2})_n},$$ is increasing as a 
function of $0\le n\le L$ and therefore 
$$ N^{\frac{d+1}{d}} \mathbb{E} [\mbox{\rm wce} (X_N;\mathbb{H}^s(\mathbb{S}^d))^2]\gtrsim \sum_{n=0}^{L}\frac{\left(\frac{d+1}{2}\right)_n \left(\frac{1}{2}\right)_n }{ \left(\frac{d}{2}+1\right)_n }\frac{1}{n!}.$$

Finally, this last series diverges when $L\to \infty$ by Gauss test taking $a=(d+1)/2,$ $b=1/2$ and $c=a+b$, because
$$\frac{\frac{(a)_n (b)_n}{(c)_n n!}}{\frac{(a)_{n+1} (b)_{n+1}}{(c)_{n+1} (n+1)!}}=\frac{(c+n)(n+1)}{(a+n)(b+n)}=1+\frac{1}{n}+\frac{C_n}{n^2},$$
with $C_n$ a bounded sequence.
\qed


\begin{proof}(Corollary \ref{DeMa_cor_as})
Recall that a measure of the uniformity of a finite set of points is the $L^\infty$ spherical cap discrepancy
$\mathbb{D}_\infty(X_N)$
where $X_N=\{x_1, \dots, x_N\}\subset \mathbb{S}^d.$ It was proved in \cite{BeMaOr} that for every $M>0$ there exist $C_M>0$ such that
$$\mathbb P\left( \mathbb{D}_\infty(X_N)\le C_M \frac{\log N}{N^{\frac{d+1}{2d}}}  \right)\ge 1-\frac{1}{N^M},$$
where $X_N$ is an $N$-point set drawn from the harmonic ensemble.

Now it follows from the result above, formula (\ref{DeMa_equality_extended}) and 
$\mathbb{D}_2 (X_N)\lesssim \mathbb{D}_\infty(X_N)$ that there exists (another)
$C_M >0$ such that
$$\mathbb P\left( \mbox{\rm wce} (X_N;\mathbb{H}^{\frac{d+1}{2}}(\mathbb{S}^d)) \le C_M \frac{\log N}{N^{\frac{d+1}{2d}}}  \right)\ge 1-\frac{1}{N^M}.$$
To get (\ref{DeMa_bounds_prob}) it is enough to apply the interpolation result from \cite[Lemma 23]{BrSaSlWo}, from which we get for $d/2<s<(d+1)/2$ a constant $C_{d,s}>0$ such that 
$\mbox{\rm wce} (X_N;\mathbb{H}^{s}(\mathbb{S}^d))^{\frac{d+1}{2s}}\le C_{d,s}\mbox{\rm wce} (X_N;\mathbb{H}^{\frac{d+1}{2}}(\mathbb{S}^d))$
if $\mbox{\rm wce} (X_N;\mathbb{H}^{\frac{d+1}{2}}(\mathbb{S}^d))\le 1$, see \cite[Section 1.5]{Be}. 

Finally, if we take (for example) $M=2$ in (\ref{DeMa_bounds_prob}) we get
$$\sum_{N} \mathbb P\left( \mbox{\rm wce} (X_N;\mathbb{H}^{s}(\mathbb{S}^d))  > C_{d,s} \frac{(\log N)^{\frac{2s}{d+1}}}{N^{\frac{s}{d}}}  \right)<\infty,$$ 
and from Borel-Cantelli lemma
$$P\left(  \limsup_{N\to +\infty} \left\{ \mbox{\rm wce} (X_N;\mathbb{H}^{s}(\mathbb{S}^d))  > C_{d,s} \frac{(\log N)^{\frac{2s}{d+1}}}{N^{\frac{s}{d}}} \right\} \right)=0.$$
\end{proof}

\subsection{The spherical ensemble (Theorem \ref{DeMa_wce_spherical})}
\begin{proof}(Theorem \ref{DeMa_wce_spherical})
The first part is due to \cite{Hi}. We include it for the sake of completeness. Let $s\in(1,2)$. Then the worst-case error is given by \eqref{DeMa_s12}.
Taking expectations and using \eqref{DeMa_energy_sph} with $s'=2-2s$,
\[
\begin{aligned}
N^s\mathbb{E}\left[\mbox{wce}(X_N,\mathbb{H}^{s}(\mathbb{S}^2))^2\right]
&=
-\frac{N^s}{N^2}\left(\mathbb{E}\left[E_{2-2s}\right]-V_{2-2s}(\mathbb{S}^2)N^2\right)\\
&=
\frac{2^{2s}\Gamma(s)}{4}\frac{\Gamma(N)}{\Gamma(N+s)}N^s\\
&\xrightarrow[N\to\infty]{}\frac{2^{2s}\Gamma(s)}{4},
\end{aligned}
\]
since $\frac{\Gamma(N)}{\Gamma(N+s)}\sim N^{-s}$ by the asymptotic property of the gamma function. Then $N^s\mathbb{E}\left[\mbox{wce}(X_N,\mathbb{H}^{s}(\mathbb{S}^2))^2\right]$ is bounded and $(X_N)$ is a sequence of $s$-QMC designs on average for $s\in(1,2)$.

Now we show that $s^{*}=2$. Let $s\in(2,3)$. The expression for the worst-case error is \eqref{DeMa_others} with $M=1$:
\[\mbox{wce}(X_N;\mathbb{H}^s(\mathbb{S}^2))^2=\frac{1}{N^2} \left[\sum_{j,i=1}^N \mathcal{Q}_1(x_j\cdot x_i)+E_{2-2s}-V_{2-2s}(\mathbb{S}^2)N^2\right],
\]
with
\[\mathcal{Q}_1(x_j\cdot x_i)=-6\alpha_1^{(s)}x_j\cdot x_i,\]
where we have used that $P_1^{(2)}(x)=x$.
Then
\begin{equation}\label{DeMa_s23}
\begin{aligned}
\mbox{wce}(X_N;\mathbb{H}^s(\mathbb{S}^2))^2&= \frac{1}{N^2}\left[-6\alpha_1^{(s)}\sum_{j,i=1}^N x_j\cdot x_i+E_{2-2s}-V_{2-2s}(\mathbb{S}^2)N^2\right]\\
&= \frac{1}{N^2}\left[-6\alpha_1^{(s)}\sum_{j,i=1}^N \left(1-\frac{|x_j-x_i|^2}{2}\right)+E_{2-2s}-V_{2-2s}(\mathbb{S}^2)N^2\right]\\
&= \frac{1}{N^2}\left[3\alpha_1^{(s)} \left(E_{-2}-2N^2\right)+E_{2-2s}-V_{2-2s}(\mathbb{S}^2)N^2\right]\\
\end{aligned}
\end{equation}
Taking expectations and using \eqref{DeMa_energy_sph} with $s=-2$  and  $s'=2-2s$,
\[
\begin{aligned}
N^2\mathbb{E}\left[\mbox{wce}(X_N;\mathbb{H}^s(\mathbb{S}^2))^2\right]&=\left[3\alpha_1^{(s)} \left(\mathbb{E}\left[E_{-2}\right]-2N^2\right)+\mathbb{E}\left[E_{2-2s}\right]-V_{2-2s}(\mathbb{S}^2)N^2\right]\\
&=-12\alpha_1^{(s)}\frac{N}{N+1}-\frac{2^{2s}\Gamma(s)}{4}\frac{\Gamma(N)}{\Gamma(N+s)}N^2\\
&\xrightarrow[N\to\infty]{}-12\alpha_1^{(s)}>0,
\end{aligned}
\]
since $\frac{\Gamma(N)}{\Gamma(N+s)}N^2\sim N^{2-s}$. Hence, $N^2\mathbb{E}\left[\mbox{wce}(X_N;\mathbb{H}^s(\mathbb{S}^2))^2\right]$ is bounded below by a positive constant.
\end{proof}

\subsection{Elliptic polynomials (Theorem \ref{DeMa_wce_theorem})}
\begin{proof}(Theorem \ref{DeMa_wce_theorem})
Let $s\in(2,3)$. We have already seen in \eqref{DeMa_s23} that the expression for the worst-case error is
\[
\mbox{wce}(X_N;\mathbb{H}^s(\mathbb{S}^2))^2=  \frac{1}{N^2}\left[3\alpha_1^{(s)} \left(E_{-2}-2N^2\right)+E_{2-2s}-V_{2-2s}(\mathbb{S}^2)N^2\right].
\]
Taking expectations and using \eqref{DeMa_cas-2} for $\mathbb{E}[E_{-2}]$ and \eqref{DeMa_weakasymptotics} with $s'=2-2s$,
\[
\begin{aligned}
N^s\mathbb{E}\left[\mbox{wce}(X_N;\mathbb{H}^s(\mathbb{S}^2))^2\right]&=\frac{N^s}{N^2}\left[3\alpha_1^{(s)} \left(\mathbb{E}\left[E_{-2}\right]-2N^2\right)+\mathbb{E}\left[E_{2-2s}\right]-V_{2-2s}(\mathbb{S}^2)N^2\right]\\
&=N^{s-2}\bigg[3\alpha_1^{(s)} \left(-8\zeta(3)\frac{1}{N}+o\left(\frac{1}{N}\right)\right)\\
&+C(2-2s)N^{2-s}+o(N^{2-s})\bigg]\\
&=3\alpha_1^{(s)} \left(-8\zeta(3)N^{s-3}+o\left(N^{s-3}\right)\right)+C(2-2s)+o(1)\\
&\xrightarrow[N\to\infty]{} C(2-2s).\\
\end{aligned}
\]
Then $N^s\mathbb{E}\left[wce(Q[X_L];\mathbb{H}^s(\mathbb{S}^2))^2\right]$ is bounded for $s\in (2,3)$. For $1<s\leq 2$, the result holds automatically from Proposition \ref{DeMa_lowranQMC}.

Now we see that the strength is $s^{*}=3$. Let $s\in (3,4)$. By \eqref{DeMa_others} with $M=2$, the square of the worst-case error is
\[
\mbox{wce}(X_N;\mathbb{H}^s(\mathbb{S}^2))^2
=\frac{1}{N^2} \left[\sum_{j,i=1}^N \mathcal{Q}_2(x_j\cdot x_i)-(E_{2-2s}-V_{2-2s}(\mathbb{S}^2)N^2)\right],
\]
where, using that $P_2^{(2)}(x)=\frac{1}{2}(3x^2-1)$,
\[
\begin{aligned}
\mathcal{Q}_2(x_j\cdot x_i)&=-5\alpha_2^{(s)}[3(x_j\cdot x_i)^2-1]\\
&=-5\alpha_2^{(s)}\left[2-3|x_j- x_i |^2+\frac{3}{4}|x_j- x_i |^4\right]
.
\end{aligned}\]
Then
\[
\begin{aligned}
\mbox{wce}(X_N;\mathbb{H}^s(\mathbb{S}^2))^2
=&\frac{1}{N^2} \Big[-5\alpha_2^{(s)}\Big(2N^2-3E_{-2}+\frac{3}{4}E_{-4}\Big)\\
&-(E_{2-2s}-V_{2-2s}(\mathbb{S}^2)N^2)\Big].
\end{aligned}
\]

Taking expectations and using \eqref{DeMa_cas-2} for $\mathbb{E}[E_{-2}]$ and \eqref{DeMa_weakasymptotics} with $s'=2-2s$,
\[
\begin{aligned}
\mathbb{E}\left[\mbox{wce}(X_N;\mathbb{H}^s(\mathbb{S}^2))^2\right]
=&\frac{1}{N^2} \Big[-5\alpha_2^{(s)}\Big(2N^2-3\mathbb{E}[E_{-2}]+\frac{3}{4}\mathbb{E}[E_{-4}]\Big)\\
&-(\mathbb{E}[E_{2-2s}]-V_{2-2s}(\mathbb{S}^2)N^2)\Big]\\
=&\frac{1}{N^2} \Bigg\{-5\alpha_2^{(s)}\left[2N^2-3\left(2N^2-8\zeta(3)\frac{1}{N}+o\left(\frac{1}{N}\right)\right)\right.\\
&\left.+\frac{3}{4}\left(\frac{32}{6}N^2+64\zeta(3)\frac{1}{N}+o\left(\frac{1}{N}\right)\right)\right]\\
&-\left(C(2-2s)N^{2-s}+o\left(N^{2-s}\right)\right)\Bigg\}\\
=&\frac{1}{N^2} \Bigg\{-5\alpha_2^{(s)}\left[72\zeta(3)\frac{1}{N}+o\left(\frac{1}{N}\right)\right]-C(2-2s)N^{2-s}\Bigg\}.
\end{aligned}
\]
Therefore,
\[\begin{aligned}
N^3\mathbb{E}\left[\mbox{wce}(X_N;\mathbb{H}^s(\mathbb{S}^2))^2\right]&=-5\alpha_2^{(s)}\left[72\zeta(3)+o\left(1\right)\right]-C(2-2s)N^{3-s}\\
&\xrightarrow[N\to\infty]{}-360\zeta(3)\alpha_2^{(s)}>0.
\end{aligned}
\]
Then there exists a constant $B>0$ such that
\[N^3\mathbb{E}\left[\mbox{wce}(X_N;\mathbb{H}^s(\mathbb{S}^2))^2\right]\geq B.\]
\end{proof}


\end{document}